\documentclass[12pt]{amsart}
\usepackage{geometry, amssymb}
\usepackage[all]{xy} 
\title{$K_1$ of some noncommutative $p$-adic group rings}
\author{Mahesh Kakde}
\date{} 
\newtheorem{theorem}{Theorem}
\newtheorem{lemma}[theorem]{Lemma}
\newtheorem*{definition}{Definition}
\newtheorem{corollary}{Corollary}

\newcommand{\Z}{\mathbb{Z}}
\newcommand{\zp}{\mathbb{Z}_p}
\newcommand{\C}{\mathcal{C}}

\begin{document}

\maketitle
\tableofcontents
Let $G$ be a finite $p$-group, for any prime $p$. In this short note we will describe $K_1(\mathbb{Z}_p[G])$ modulo its $p$-power torsion in terms of abelian subquotients of $G$. Such a description has application in noncommutative Iwasawa theory due to a strategy proposed by D. Burns and K. Kato with a modification due to Hara \cite{Hara:2009}. We will however, not say anything about the strategy and about noncommutative Iwasawa theory here. We must mention that such description of the Whitehead group of group rings or completed group rings in terms of abelian subquotients was first given by K. Kato \cite{Kato:2005} for completed group rings of open subgroups of $\mathbb{Z}_p\rtimes \mathbb{Z}_p^{\times}$. The proof was rather tedious and long. In an unpublished article \cite{Kato:2006} K. Kato himself found a much more elegant approach using integral logarithm of R. Oliver and M. Taylor. He used this approach to describe the Whitehead groups of completed group rings of $p$-adic Heisenberg group and also proved the ``main conjecture" for Galois extension of totally real number fields whose Galois group is a quotient of the $p$-adic Heisenberg group. The author generalised this approach in his PhD thesis \cite{Kakde:2008} to prove main conjucture for pro-$p$ groups of  ``special type". These groups include many interesting examples such as $\mathbb{Z}_p^r\rtimes\mathbb{Z}_p$ where the action of $\mathbb{Z}_p$ on $\mathbb{Z}_p^r$ is diagonal. T. Hara \cite{Hara:2009} used the approach of K. Kato and some inductive arguments to prove the main conjucture for Galois extensions of totally real number fields whose Galois group is of the form $H\times\Gamma$, where $H$ is any finite group of exponent $p$. In this note we will use K. Kato's approach using the integral logarithm for describing $K_1(\zp[G])$ modulo its $p$-power torsion for a finite $p$-group $G$. 

The author thanks Professor John Coates very much for his constant encouragement. Thanks are due to the hospitality of Newton Institute, Cambridge where this work was done.

\section{The additive side}
Let $G$ be a finite group and let $Conj(G)$ be the set of conjugacy classes of $G$. Let $\zp[Conj(G)]$ be the free $\zp$ module generated by $Conj(G)$. We first describe $\zp[Conj(G)]$ in terms of abelian subquotients of $G$. If $H$ is any subgroup of $G$, then $C(G,H)$ denotes a set of left coset representatives of $H$ in $G$. Unless stated otherwise we will use the notation $[g]$ to denote the conjugacy class of $g$ (the group will be clear from the context or stated explicitly) and use $\bar{g}$ to denote the image of $g$ in abelianisation (again the group will be clear from the context or stated explicitly). Let $N_GH$ denote the normaliser of $H$ in $G$ and let $W_GH:=N_GH/H$ be the Wyel group. The group $N_GH$ acts by conjugation on $\zp[Conj(H)]$ and $\zp[H^{ab}]$. The action of $H \leq N_GH$ on $\zp[H^{ab}]$ is trivial and hence we get an action of $W_GH$ on $\zp[H^{ab}]$. For any $H \leq G$ we have the following trace map on $\zp[H^{ab}]$ given by 
\[
x \mapsto \sum_{g \in W_GH} gxg^{-1}.
\]
Let $T_H$ be the image of this map. Thus $T_H$ is an ideal in $\zp[H^{ab}]^{W_GH}$. 

\subsection{Maps relating $\zp[Conj(G)]$ to abelian subquotients} Let $H$ be a subgroup of $G$. Then we define a map $t^G_H:\zp[Conj(G)] \rightarrow \zp[Conj(H)]$ by
\[
t^G_H([g]) = \sum_{x \in C(G,H)}\{[x^{-1}gx] : x^{-1}gx \in H\}
\]
This induces a well defined map 
\[
\beta^G_H : \zp[Conj(G)] \rightarrow \zp[H^{ab}]
\]
Let $\beta^G$ be the map 
\[
\beta^G := (\beta^G_H)_{H\leq G} : \zp[(G)] \rightarrow \prod_{H\leq G} \zp[H^{ab}].
\]
Let $\C$ be the set of cyclic subgroups of $G$. Let $\beta^G_{\C}$ be the map
\[
\beta^G_{\C} := (\beta^G_H)_{H\in \C} = \zp[Conj(G)] \rightarrow \prod_{H \in \C} \zp[H].
\]
The maps $\beta^G$ and $\beta^G_{\C}$ fit in the following commutative diagram 
\[
\xymatrix{ \zp[Conj(G)] \ar[r]^{\beta^G} \ar[dr]_{\beta^G_{\C}} & \prod_{H \leq G} \zp[H^{ab}] \ar[d]^{proj} \\ & \prod_{H \in \C} \zp[H] }
\]

\subsection{The image of $\beta^G_{\C}$} 

\begin{definition} Let $\Phi_{\C}$ be the $\zp$-submodule of $\prod_{H\in \C}\zp[H]$ consisting of all tuples $(a_H)_{H \in \C}$ satisfying \\
\noindent A1) For any $H \leq H^{\prime} \leq G$, with $H, H^{\prime} \in \C$, we have that 
\[
\beta^{H^{\prime}}_H (a_{H^{\prime}}) =a_H.
\]

\noindent A2) For any $g \in G$, we want $(a_H)_{H \in \C}$ to be fixed by $g$ under the conjugation action. In particular for any $H \in \C$, we have $a_H \in \zp[H]^{N_GH}$.  \\

\noindent A3) For any $H \in \C$, we want $a_H \in T_H$. 
\end{definition} 

\begin{lemma} The image of $\beta^G_{\C}$ is contained in $\Phi_{\C}$. 
\label{lemma1}
\end{lemma}
\noindent{\bf Proof:} It is enough to show that $\beta^G_{\C}([g]) \in \Phi_{\C}$ for any $g\in G$, i.e. that it satisfies A1), A2) and A3).  \\

\noindent A1): $\beta^G_{\C}([g])$ satisfies A1) because for any $H \leq H^{\prime} \leq G$, with $H, H^{\prime} \in \C$, we have the following commutative diagram 
\[
\xymatrix{ \zp[Conj(G)] \ar[r]^-{\beta^G_{H^{\prime}}} \ar[rd]_{\beta^G_H} & \zp[H^{\prime}] \ar[d]^{\beta^{H^{\prime}}_H} \\ & \zp[H] }
\]
\noindent A2): Let $g, g_1 \in G$. We must show that $g_1\beta^G_H([g])g_1^{-1} = \beta^G_{g_1Hg_1^{-1}}([g])$ for any $H \in \C$. We in fact show that this holds for any $H \leq G$. 
\begin{align*}
g_1\beta^G_H([g])g_1^{-1} &= g_1(\sum_{x \in C(G,H)} \{ \overline{[x^{-1}gx]} : x^{-1}gx \in H\}) g_1^{-1} \\
&= \sum_{x\in C(G,H)} \{ \overline{[(g_1x^{-1}g_1^{-1})(g_1gg_1^{-1})(g_1xg_1^{-1})]} : x^{-1}gx \in H\} \\
&= \sum_{x_1 \in C(G,g_1Hg_1^{-1})} \{ \overline{[x_1^{-1} (g_1gg_1^{-1})x_1]} : x_1^{-1}g_1gg_1^{-1}x_1 \in g_1Hg_1^{-1}\} \\
&= \beta^G_{g_1Hg_1^{-1}}([g_1gg_1^{-1}]) \\
&=\beta^G_{g_1Hg_1^{-1}}([g]).
\end{align*}

\noindent A3): For A3) we must show that $\beta^G_H([g]) \in T_H$ for any $H \in \C$. Again we show this for any $H \leq G$. First note that
\[
\beta^G_H([g]) = \beta^{N_GH}_H(\beta^G_{N_GH}([g])) 
\]
Hence it is enough to show that $\beta^{N_GH}_H([g]) \in T_H$ for any $g\in N_GH$. But $\beta^{N_GH}_H([g])$ is non-zero if and only is $g \in H$, and when $g\in H$, we have
\begin{align*}
\beta^{N_GH}_H([g]) &= \sum_{x \in C(N_GH, H)} x^{-1}gx \\
&= \sum_{x \in W_GH} x^{-1}gx \in T_H
\end{align*}
This finishes proof of the lemma. \qed \\

\noindent We now define a left inverse of $\beta^G_{\C}$. Define 
\[
\tau^G_{\C} : \prod_{H \in \C} \zp[H] \rightarrow \mathbb{Q}_p[Conj(G)],
\]
by defining 
\[
\tau^G_{\C, H} : \zp[H] \rightarrow \mathbb{Q}_p[Conj(G)],
\]
for each $H \in \C$ and putting $\tau^G_{\C} = \sum_{H \in \C} \tau^G_{\C, H}$. Define $\tau^G_{\C,H}$ by 
\[
\tau^G_{\C,H}(h) = \left\{ 
\begin{array}{l l}
\frac{1}{[G:N_GH]|W_GH|}[h] & \quad \mbox{if $h$ is a generator of $H$} \\
0 & \quad \mbox{if not} \\
\end{array} \right.
\]
For any non-trivial cyclic group $H$ of $p$-power order, let $\omega_H$ denote any nontrivial character of $H$ of order $p$. We fix such a character $\omega_H$ for each cyclic subgroup of a $p$-power order. Let $\eta_H : \zp[H] \rightarrow \zp[H]$ be the map defined by 
\[
\eta_H(h) = h - \frac{1}{p}\sum_{k=0}^{p-1}\omega_H^k(h)h = \left \{ \begin{array}{ll} 
h &  \quad \mbox{if $h$ generates $H$} \\
0  & \quad \mbox{otherwise} 
\end{array} \right. 
\]
And for the trivial group $H$ we take $\eta_H$ to be the identity function. Then we may define $\tau^G_{\C,H}$ as 
\[
\tau^G_{\C,H}(h)= \frac{1}{[G:N_GH]|W_GH|}[\eta_H(h)].
\]

\begin{lemma} $\tau^G_{\C} \circ \beta^G_{\C}$ is identity on $\zp[Conj(G)]$. In particular, $\beta^G_{\C}$ is injective.
\end{lemma}
\noindent{\bf Proof:} We show that $\tau^G_{\C}(\beta^G_{\C}([g])) = [g]$, for any $g \in G$. Let $H$ be the cyclic subgroup of $G$ generated by $g$. Let $\C_H$ be the set of conjugates of $H$ in $G$. Then we note that 
\begin{align*}
\tau^G_{\C}(\beta^G_{\C}([g])) &= \sum_{P \in \C} \tau^G_{\C, P}(\beta^G_P([g])) \\
&= \sum_{P \in \C_H} \tau^G_{\C,P}(\beta^G_{P}([g])) \\
&=\sum_{x \in C(G, N_GH)} \tau^G_{\C, xHx^{-1}}(\beta^G_{xHx^{-1}}([g])) \\
&=\sum_{x \in C(G, N_GH)} \tau^G_{\C, xHx^{-1}}(\beta^G_{xHx^{-1}}([xgx^{-1}])) \\ 
&=\sum_{x \in C(G,N_GH)} \tau^G_{\C, xHx^{-1}}(\sum_{y\in W_GH} \overline{[y^{-1}xgx^{-1}y]}) \\ 
&=\sum_{x \in C(G, N_GH)} |W_GH|\tau^G_{\C, xHx^{-1}}(\overline{[xgx^{-1}]}) \\
&=\sum_{x \in C(G, N_GH)} \frac{1}{[G:N_GH]} [xgx^{-1}] = [g] . 
\end{align*} \qed

\begin{lemma} $\tau^G_{\C}$ restricted to $\Phi_{\C}$ is injective and its image lies in $\zp[Conj(G)]$. 
\end{lemma}
\noindent{\bf Proof:} Let $(a_H)_{H\in \C} \in \Phi_{\C}$ be such that $\tau^G_{\C}((a_H)) = \sum_{H\in \C} \tau^G_{\C, H} (a_H) = 0$. We claim that $\tau^G_{\C, H}(a_H)$ is zero for every $H$. This follows from two simple observations: firstly, $\tau^G_{\C, H}(a_H)$ and $\tau^G_{\C, H^{\prime}}(a_{H^{\prime}})$ cannot cancel each other unless $H$ and $H^{\prime}$ are conjugates. But if $H$ and $H^{\prime}$ are conjugates of each other, then by A2)
\[
\tau^G_{\C, H}(a_H) = \tau^G_{\C, H^{\prime}}(a_{H^{\prime}}).
\]
Hence $\tau^G_{\C, H}(a_H)=0$ for any $H \in \C$. This implies that the coefficients of generators of $H$ in $a_H$ are 0. Now note that if $H \leq H^{\prime}$ are two cyclic subgroups of $G$ and if $a_{H^{\prime}} = \sum_{h\in H^{\prime}} a_{H^{\prime}, h}h$, then 
\[
\beta^{H^{\prime}}_H(a_{H^{\prime}}) = [H^{\prime}:H]\sum_{h \in H}a_{H^{\prime},h}h. 
\]
By A1) this must be equal to $a_H$. Hence none of the generators of $H$ appear in $a_{H^{\prime}}$. This is true for any subgroup of $H^{\prime}$. Hence $a_{H^{\prime}} =0$. Now we prove the second claim of the lemma. Let $(a_H)_{H\in \C} \in \Phi_{\C}$. Then $a_H \in T_H$, for all $H\in \C$ by A3). Let $a_H= \sum_{x \in W_GH} xb_Hx^{-1}$ for some $b_H \in \Z_p[H]$. Then
\begin{align*} 
\tau^G_{\C}((a_H)) &= \sum_{H \in \C}\tau^G_{\C, H}(a_H) \\
&=\sum_{H \in \C/G} [G:N_GH] \tau^G_{\C,H}(\sum_{x\in W_GH}xb_Hx^{-1}) \\
&=\sum_{H\in \C/G}[G:N_GH]|W_GH| \tau^G_{\C, H}(b_H) \in \zp[Conj(G)].
\end{align*}
\qed

\begin{theorem} $\beta^G_{\C}$ induces an isomorphism between $\zp[Conj(G)]$ and $\Phi_{\C}$. 
\end{theorem}
\noindent{\bf Proof:} Let $\tau = \tau^G_{\C}|_{\Phi_{\C}}$. Then we know that $\tau \circ \beta^G_{\C} $ is identity on $\zp[Conj(G)]$. We claim that $\beta^G_{\C} \circ \tau$ is identity on $\Phi_{\C}$.  Let $(a_H)_{H\in \C} \in \Phi_{\C}$, then $\tau(\beta^G_{\C}(\tau((a_H))))= \tau((a_H))$. Hence by the previous lemma $\beta^G_H(\tau((a_H))) = (a_H)$. \qed 

\subsection{The image of $\beta^G$} Now lets assume that $G$ is a $p$-group. Let $H \leq H_1$ be two subgroups of $G$. Assume that $[H_1, H_1] \leq H$, then $H/[H_1,H_1] \leq H_1^{ab}$ and  we have two naturally defined maps
\[
\mbox{trace} = tr_{H,H_1} : \zp[H_1^{ab}] \rightarrow \zp[H/[H_1,H_1]],
\]
and the natural surjection
\[
\pi_{H,H_1} = \zp[H^{ab}] \rightarrow \zp[H/[H_1,H_1]].
\]

\begin{definition} Let $\Phi^G$ be the $\zp$-submodule of $\prod_{H\leq G} \zp[H^{ab}]$ consisting of all tuples $(a_H)_{H \leq G}$, such that \\
A1) For any two subgroups $H \leq H_1$ of $G$ such that $[H_1, H_1] \leq H$, we have 
\[
tr_{H, H_1}(a_{H_1}) = \pi_{H,H_1}(a_H).
\]
\noindent A2) We want $(a_H)_{H \leq G}$ to be fixed by all $g \in G$. In particular, $a_H \in \zp[H^{ab}]^{W_GH}$ for all $H \leq G$. \\
\noindent A3) We want $a_H \in T_H$ for all $H \in \C$. 
\end{definition}

\begin{lemma} The image of $\beta^G$ is contained in $\Phi^G$. 
\end{lemma}
\noindent{\bf Proof:} We have already shown that image of $\beta$ satisfies A2) and A3) in the proof of lemma \ref{lemma1}. In fact, we showed that A2) and A3) is satisfied for every $H \leq G$. We show that it satisfies A1). First note that $[H_1,H_1] \leq H$ implies that $H$ is a normal subgroup of $H_1$ (since if $h\in H$ and $x \in H_1$, then $xhx^{-1}h^{-1} \in [H_1,H_1] \leq H$. Hence $xhx^{-1} \in Hh=H$). Now we must show that the following diagram commutes.
\[
\xymatrix{\zp[Conj(G)] \ar[r]^-{\beta^G_{H_1}} \ar[d]_-{\beta^G_H} & \zp[H_1^{ab}] \ar[d]^-{tr_{H,H_1}} \\ \zp[H^{ab}] \ar[r]_-{\pi_{H,H_1}} & \zp[H/[H_1,H_1]]}
\]
Let $h \in H_1$, then $tr_{H,H_1}(\bar{h})$ is 0 unless $h \in H$ in which case it is 
\[
tr_{H,H_1}(\bar{h})=[H_1^{ab}:(H/[H_1,H_1])]\bar{h} = [H_1:H]\bar{h} \in \zp[H/[H_1,H_1]].
\]
On the other hand $\pi_{H,H_1}(\beta^{H_1}_H([h]))$ is 0 unless $h\in H$ in which case it is
\[
\pi_{H,H_1}(\beta^{H_1}_H([h]))= [H_1:H]\bar{h}\in \zp[H/[H_1,H_1]].
\]
Hence for any $g\in G$, we have 
\begin{align*}
tr_{H_1,H}(\beta^G_{H_1}([g]))&= [H_1:H] \sum_{x \in C(G,H_1)} \{ \overline{[x^{-1}gx]} : x^{-1}gx \in H\} \\
&= \pi_{H_1,H}(\beta^G_H([g]))
\end{align*} 
\qed

\begin{lemma} The projection $proj : \prod_{H \leq G} \zp[H^{ab}] \rightarrow \prod_{H \in \C} \zp[H]$ maps $\Phi^G$ isomorphically onto $\Phi_{\C}$
\end{lemma}
\noindent{\bf Proof:} Surjectivity follows from $proj(\beta^G((\beta^G_{\C})^{-1}((a_H)))) = (a_H)$. We now prove injectivity of $proj : \Phi^G \rightarrow \Phi_{\C}$. Let $(a_H)_{H\leq G}$ be such that $proj((a_H)) =0 $ in $\Phi_{\C}$. We will use induction on order of $H$. Let $H \leq G$ be a non-cyclic subgroup of $G$. We must show that $a_H = 0$. Let $a_H = \sum_{h \in H^{ab}} a_hh$. Let $h_0 \in H^{ab}$ be such that $a_{h_0} \neq 0$. Let $\tilde{h_0}$ be any lift of $h_0$ to $H$. Let $P$ be a maximal subgroup of $H$ containing $\tilde{h_0}$. Then $P$ is a normal subgroup of $H$ of index $p$ (Since $H$ is a non-cyclic $p$-group). By the induction hypothesis $a_P = 0$. By A1) $tr_{P, H}(a_H) = 0$. But the co-efficient of $h_0 \in P/[H,H]$ in $tr_{P,H}(a_H)$ is $pa_{h_0} \neq 0$. This contradicts $tr_{P,H}(a_H) = 0$. \qed  \\

\noindent{\bf Remark:}  We use the hypothesis that $G$ is a $p$-group for the first time in the above lemma.

\begin{theorem} Let $G$ be a finite $p$ group. Then $\beta^G$ induces an isomorphism between $\zp[Conj(G)]$ and $\Phi^G$.
\end{theorem}
\noindent{\bf Proof:} This is a straightforward corollary of the previous lemma. \qed \\

\noindent We now explicitly construct the inverse $q$ of $proj: \Phi^G \rightarrow \Phi_{\C}$. First define $q_H : \Phi_{\C} \rightarrow \zp[H^{ab}]$ for each $H \leq G$ and then put $q = (q_H)_{H\leq G}$. For any $H\leq G$, define $q_H$ by 
\[
q_H((a_P)) = \sum_{P \leq H \mbox{ and } H\in \C}\frac{1}{[H:N_HP]|W_HP|} \overline{\eta_P(a_p)},
\]
where $\overline{\eta_P(a_P)}$ denote the image of $\eta_P(a_P)$ under the natural map $P \rightarrow H^{ab}$. 

\begin{lemma} Image of $q$ is contained in $\Phi^G$ and $q$ is the inverse of $proj$. 
\end{lemma}
\noindent{\bf Proof:} Let $\beta^G_{\C}([g]) = (a_P)_{P \in \C}$. Both the claims in lemma follow if we show that $q_H((a_P)) = \beta^G_H([g])$ for all $H \leq G$. Note that $\eta_P(a_P)$ is non-zero if and only if $P$ is a conjugate of the cyclic group generated by $g$. Let $P_g$ be the cyclic group generated by $g$ and let $\C_g$ be the set of all conjugates of $P_g$. Let $\C_{g,H}$ consist of all $P \in \C_g$ which are subgroups of $H$.Then 
\begin{align*}
q_H((a_P)) &= \sum_{P\in\C_{g,H}} \frac{1}{[H:N_HP]|W_HP|}\overline{\eta_P(a_P)}\\
&=\sum_{P\in\C_{g,H}} \frac{1}{[H:N_HP]|W_HP|}\overline{\eta_P(\beta^G_P([g]))} \\
&=\sum_{P\in\C_{g,H}} \frac{1}{[H:N_HP]|W_HP|}\overline{\beta^G_P([g])} \quad (\mbox{since } \eta_P(\beta^G_P([g]))=\beta^G_P([g]))\\
&=\sum_{P\in\C_{g,H}/H}\frac{1}{|W_HP|}\overline{\beta^G_P([g])} \\
&=\sum_{P\in\C_{g,H}/H}\frac{1}{|W_HP|}\overline{\beta^H_P(\beta^G_H([g]))} \\
&=\sum_{P\in\C_{g,H}/H}\frac{1}{|W_HP|} \beta^H_P(\sum_{x\in C(G,H)}\{\overline{[x^{-1}gx]}|x^{-1}gx \in H\}) \\
&=\sum_{P\in\C_{g,H}/H}\sum_{x\in C(G,H)} \{\overline{[x^{-1}gx]}| x^{-1}gx \in P\} \\
&=\sum_{x\in C(G,H)} \{\overline{[x^{-1}gx]} | x^{-1}gx \in H\} \\
&=\beta^G_H([g]).
\end{align*}
\qed

For any $(a_P)_{P\in \C} \in \Phi_{\C}$ and any $H \in \C$, define 
\[
v_H : \Phi_{\C} \rightarrow \zp[H],
\]
by 
\[
v_H((a_P)) = \sum_{P^p \leq H} \frac{[P:P^p]}{[H:P^p]}ver_{P,P^p}(\eta_P(a_P)) \in \zp[H],
\]
where $ver_{P,P^p}:P \rightarrow P^p \hookrightarrow H$ is  just the $p$-power map. Put 
\[
v = (v_H)_{H\in \C} : \Phi_{\C} \rightarrow \prod_{H\in \C} \zp[H].
\]

\begin{lemma} We have the following commutative diagram
\[
\xymatrix{ \zp[Conj(G)] \ar[d]_-{\varphi} \ar[r]^-{\beta^G_{\C}} &  \Phi_{\C} \ar[d]^-{v} \\
\zp[Conj(G)] \ar[r]_-{\beta^G_{\C}} & \Phi_{\C}}
\]
\end{lemma}
\noindent{\bf Proof:} We will show that for any $H\in \C$, we have $\beta^G_{\C,H}([g^p]) = v_H(\beta^G_{\C}([g]))$. It is best to consider the following two cases: \\
Case 1: $g=1$. \quad Then $\beta^G_{\C, H}([1])=\frac{|G|}{|H|}$. On the other hand
\[
v_{H}(\beta^G_{\C}([1])) = \frac{1}{[H:\{1\}]}(|G|) = \frac{|G|}{|H|}.
\]
Case 2: $g \neq 1$. Let $P_g$ be the cyclic group generated by $g$. Let $\C_g$ be all the conjugates of $P_g$ in $G$. Then
\begin{align*} 
v_H(\beta^G_{\C}([g]))&=\sum_{P\in \C_g, P^p \leq H}\frac{p}{[H:P^p]}ver_{P,P^p}(\beta^G_{\C,P}([g])) \\
&=\sum \frac{p}{[H:p^p]} \sum_{x\in C(G,P)} \{\overline{[x^{-1}g^px]} | x^{-1}gx \in P\} \\
&=\sum \frac{1}{[H:P^p]}\sum_{x\in C(G,P^p)} \{\overline{[x^{-1}g^px]}|x^{-1}gx\in P\} \\
&=\sum \frac{1}{[H:P^p]} \sum_{x\in C(G,H)} \sum_{y \in C(H,P^p)} \{\overline{[y^{-1}x^{-1}g^pxy]}| y^{-1}x^{-1}gxy \in P\} \\
&= \sum \sum_{x\in C(G,H)} \{\overline{[x^{-1}g^px]} | x^{-1}gx \in P\} \\
&= \sum_{x \in C(G,H)} \sum_{P\in \C_g, P^p \leq H} \{\overline{[x^{-1}g^px]} | x^{-1}gx \in P\} 
\end{align*}
Now $x^{-1}gx$ lies in $P$ and $P^{\prime}$ for $P, P^{\prime} \in \C_g$ implies that $P = P^{\prime}$ and $x^{-1}gx \in P$ for some $P \in \C_g$ such that $P^p \leq H$ if and only if $x^{-1}g^px \in H$. Hence 
\[
\sum_{P \in \C_g, P^p \leq H} \{\overline{[x^{-1}g^px]} | x^{-1}gx \in P\} = \left\{ \begin{array}{l l}
\overline{[x^{-1}g^px]} & \mbox{if $x^{-1}gx \in H$}, \\
0 & \mbox{if $x^{-1}gx \notin H$}. \end{array} \right.
\]
Putting this in the equation above we obtain
\begin{align*}
v_H(\beta^G_{\C}([g]))&= \sum_{x \in C(G,H)} \{\overline{[x^{-1}g^px]}| x^{-1}gx \in H\} \\
&= \beta^G_{\C,H}([g^p]).
\end{align*} 
 \qed
 
\begin{definition} Define $v_G = q \circ v \circ proj : \Phi^G \rightarrow \Phi^G$. 
\end{definition}

\begin{corollary} The following diagram commutes
\[
\xymatrix{ \zp[Conj(G)] \ar[d]_-{\varphi} \ar[r]^-{\beta^G} & \Phi^G \ar[d]^-{v_G} \\
\zp[Conj(G)] \ar[r]_-{\beta^G} & \Phi^G}
\]
\end{corollary}
\noindent{\bf Proof:} We break the square into following two squares which we know commute
\[
\xymatrix{ \zp[Conj(G)] \ar[d]_{\varphi} \ar[rr]^{\beta^G_{\C}} & & \Phi_{\C} \ar[d]_v \ar[r] & \Phi^G \ar[d]^{v_G} \\
\zp[Conj(G)] \ar[rr]_{\beta^G_{\C}} & & \Phi_{\C} \ar[r] & \Phi^G}
\]
\qed

\noindent Explicitly, the $H^{th}$ component, $v_{G,H}$ of $v_G$ is given by 
\[
v_{G,H}((a_P)) = \sum_{P \in \C \mbox{ and } P^p \leq H} \frac{[P:P^p]}{[H:N_HP^p]|W_HP^p|} \overline{ver_{P,P^p}(\eta_P(a_P))}.
\]

\section{The multiplicative side} We have the following multiplicative analogues of $\beta^G_H$.
\[
\theta^G_H : K_1(\zp[G]) \xrightarrow{N^G_H} K_1(\zp[H]) \xrightarrow{can} \zp[H^{ab}]^{\times},
\]
 where $N^G_H$ is the norm map and $can$ is the map induced by natural surjection $H \rightarrow H^{ab}$. Let 
 \[
 \theta^G = (\theta^G_H)_{H \leq G} : K_1(\zp[G]) \rightarrow \prod_{H\leq G}\zp[H^{ab}]^{\times}.
 \]
 
\subsection{Logarithm and integral logarithm} In this subsection we recall the logarithm and integral logarithm maps of Oliver and Taylor \cite{Oliver:1988}. 

\begin{theorem} Let $J \subset \zp[G]$ be the Jacobson radical. Let $I \subset \zp[G]$ be any two sided ideal. Then the $p$-adic logarithm $Log(1+x)$ induces a unique homomorphism 
\[
log_I : K_1(\zp[G], I) \rightarrow \mathbb{Q} \otimes_{\Z}(I/[\zp[G], I]).
\]
$Ker(log_I)$ is finite for any $I$. If $I^p \subset pIJ$, then $log_I$ is an isomorphism.
\end{theorem}
This is theorem 2.8 and 2.9 in \cite{Oliver:1988}

\begin{lemma} For any $H \leq G$, we have $\beta_H \circ log = log \circ \theta_H$.
\end{lemma}
This is lemma in theorem 1.4 \cite{OliverTaylor:1988}

\begin{definition} Let $log = log_{\zp[G]}$. Define the following map
\[
L = log - \frac{\varphi}{p}log : K_1(\zp[G]) \rightarrow \mathbb{Q}_p[Conj(G)].
\]
\end{definition} 

\begin{theorem} For any finite $p$-group $G$, the image of $L$ is contained in $\zp[Conj(G)]$. The map $L$ is natural with respect to maps induced by group homomorphisms. 
\end{theorem}
This is theorem 6.2 in \cite{Oliver:1988}

\begin{theorem} Let $G$ be a finite $p$-group. Let $\epsilon = (-)^{p-1}$ and define
\[
\omega : \zp[Conj(G)] \rightarrow \langle \epsilon \rangle \times G^{ab} \qquad \mbox{by} \qquad \omega(\sum a_ig_i) = \prod(\epsilon g_i)^{a_i}.
\]
Then the following sequence is exact
\[
1 \rightarrow K_1(\zp[G])/torsion \xrightarrow{L} \zp[Conj(G)] \xrightarrow{\omega} \langle \epsilon \rangle \times G^{ab} \rightarrow 1
\]
Moreover, the torsion subgroup of $K_1(\zp[G])$ is precisely $\langle \epsilon \rangle \times \mu_{p-1} \times G^{ab} \times SK_1(\zp[G])$. Here $SK_1(\zp[G])$ is by definition the kernel of the natural map $K_1(\zp[G]) \rightarrow K_1(\mathbb{Q}_p[G])$.
\end{theorem}
This is theorem 6.6 and theorem 7.3 in \cite{Oliver:1988} (the assertion about torsion part of $K_1(\zp[G])$ is a result of C.T.C Wall and theorem 7.3 is loc. cit.).
 
\subsection{Relation between the multiplicative and additive sides} 

\begin{lemma} For any non-trivial cyclic group $P$, there is a map $\alpha_P : \Z_p[P]^{\times} \rightarrow \Z_p[P]^{\times}$ such that the diagram
\[
\xymatrix{ \Z_p[P]^{\times} \ar[r]^{log} \ar[d]_{\alpha_P} & \mathbb{Q}_p[P] \ar[d]^{p\eta_P} \\
\Z_p[P]^{\times} \ar[r]_{log} & \mathbb{Q}_p[P] }
\]
commutes.
\end{lemma}
\noindent{\bf Proof:} Since $p\eta_P(h) = ph - \sum_{k=0}^{p-1} \omega_P^k(h)$, we may define $\alpha_P$ as 
\[
\alpha_P(x) = \frac{x^p}{\prod_{k=0}^{p-1}\omega^k_P(x)}.
\]
The commutativity of the diagram in the lemma can now be verified easily since log commutes with  ring homomorphisms of $\mathbb{Z}_p[P]$ induced by homomorphisms of the group $P$.
\qed

\begin{definition} Define $\alpha_{\{1\}}$ to be the identity map. 
\end{definition}

\begin{definition} Define the map $u_{G,H} : \prod_{P \in \C} \Z_p[P]^{\times} \rightarrow \mathbb{Z}_p[H^{ab}]^{\times}$ by
\[
u_{G,H}((x_P)) \prod_{P\in \C, P^p \leq H} \overline{ver_{P,P^p}(\alpha_P(x_P))^{|P^p|}}.
\]
\end{definition} 

\begin{lemma} For any $x \in K_1(\mathbb{Z}_p[G])$ and any $H \leq G$, we have
\[
\beta^G_H(L(x)) = \frac{1}{p|H|} log\Big(\frac{\theta^G_H(x)^{p|H|}}{u_{G,H}(\theta^G(x))}\Big)
\]
\end{lemma}
\noindent{\bf Proof:} First consider $log(u_{G,H}((x_P)))$
\begin{align*} 
log(u_{G,H}((x_P))) &= \sum_{P \in \C, P^p \leq H} |P^p| ver_{P,P^p}(log(\alpha_P(x_P))) \\
&= \sum |P^p| ver_{P,P^p}([P:P^p]\eta_P(log(x_P))) \\
&= \sum |P| ver_{P,P^p} ver_{P,P^p}(\eta_P(log(x_P))) \\
&= |H| v_{G,H}((log(x_P))).
\end{align*}
Then 
\begin{align*} 
\beta^G_H(L(x)) &= \beta^G_H(log(x)) - \frac{1}{p} \beta^G_H(\phi(log(x))) \\
&= log(\theta^G_H(x)) - \frac{1}{p} v_{G,H}(\beta^G(log(x))) \\
&= log(\theta^G_H(x)) - \frac{1}{p} v_{G,H}(log(\theta^G(x))) \\
&= log(\theta^G_H(x)) - \frac{1}{p|H|} log(u_{G,H}(\theta^G(x))) \\
&= \frac{1}{p|H|} log\Big(\frac{\theta^G_H(x)^{p|H|}}{u_{G,H}(\theta^G(x))} \Big).
\end{align*}
\qed

\begin{definition} For any finite $p$-group $H$, define $J_H$ to be the kernel of the homomorphism $\Z_p[H^{ab}] \rightarrow \mathbb{F}_p$.
\end{definition}

\begin{lemma} For any $x \in K_1(\mathbb{Z}_p)$ and any $H \leq G$, we have
\[
\theta^G_H(x)^{p|H|} \equiv u_{G,H}(\theta^G(x)) \mbox{(mod $J_H$)},
\]
where $J_H = ker(\mathbb{Z}_p[H^{ab}] \rightarrow \mathbb{F}_p)$.
\label{easycongruence}
\end{lemma}
\noindent{\bf Proof:} Note that if $P \neq \{1\}$, then $\alpha_1 \equiv 1 \mbox{(mod $J_P$)}$, for any $x \in \mathbb{Z}_p[P]^{\times}$. Hence we must show that
\[
\theta^G_H(x)^{p|H|} \equiv \theta^G_{\{1\}}(x) \mbox{(mod $J_H$)},
\]
for any $H \leq G$. We must show that the following diagram commutes
\[
\xymatrix{K_1(\mathbb{F}_p[G]) \ar[d]_{norm} \ar[rdr]^{norm} & &  \\
K_1(\mathbb{F}_p[H]) \ar[rr]^{norm} \ar[d]_{\pi} &  & \mathbb{F}_p^{\times} \ar[d]^{=} \\
\mathbb{F}_p[H^{ab}]^{\times} \ar[rr]^{aug} & & \mathbb{F}_p^{\times}}.
\]
because, then $aug(\theta^G_H(x)) \equiv \theta^G_{\{1\}}(x) \mbox{(mod $p$)}$, which is what we want. To show that the diagram commutes we only have to show that the square in the lower half of the diagram commutes. We note that for any $x \in K_1(\Z_p[H])$
\[
x^{|H|} \equiv aug(x)^{|H|} \mbox{(mod $p$)},
\]
and 
\[
aug(x) = aug(\pi(x)).
\]
Hence
\begin{align*}
norm(x)^{|H|} &\equiv aug(x)^{|H|} \\
&\equiv aug(\pi(x))^{|H|} \mbox{(mod $p$)}.
\end{align*}
which show the required commutativity of the above square.
\qed

\subsection{The main theorem}

\begin{definition} Let $\tilde{K}_1(R)$ denote the group $K_1(R)/K_1(R)(p)$, where $K_1(R)(p)$ denotes the $p$-power torsion subgroup of $K_1(R)$.
\end{definition}

\begin{definition} Let denote by $\tilde{\theta}^G$ the map 
\[
\tilde{K}_1(\zp[G]) \rightarrow \prod_{H \leq G} \tilde{K}_1(\zp[H^{ab}]),
\]
induced by $\theta^G$. 
\end{definition}

\noindent The integral logarithm map $L$ is trivial on torsion and then factors through $\tilde{K}_1$. We denote the induced map from $\tilde{K}_1$ by $L$.

\begin{definition} For any subgroups $H \leq H_1 \leq G$ such that $[H_1,H_1] \leq H$ we have the norm map 
\[
nr_{H,H_1} : \tilde{K}_1(\zp[H_1^{ab}]) \rightarrow \tilde{K}_1(\zp[H^{ab}]).
\]
\end{definition}

\begin{definition} Let $\Psi^G \leq \prod_{H \leq G} \tilde{K}_1(\mathbb{Z}_p[H^{ab}])$ be the subgroup consisting of all tuples $(x_H)$ such that \\
M1) For any $H \leq H_1 \leq G$ such that $[H_1,H_1] \leq H$, we want 
\[
nr_{H,H_1}(x_{H_1}) = \pi_{H,H_1}(x_H).
\]
M2) We want all $(x_H)$ to be fixed by all $g\in G$. In particular, $x_H \in \tilde{K}_1(\Z_p[H^{ab}])^{W_GH}$, for all $H \leq G$. \\
M3) For all $P \leq G$, we want $x_P^{p|P|} \equiv u_{G,P}((x_H)) \mbox{(mod $J_P$)}$. \\
M4) For all $P \in \C$, we want $x_P^{p|P|} \equiv u_{G,P}((x_H)) \mbox{(mod $p|P|T_P$)}$.
\end{definition}

\begin{lemma} The image of $\tilde{\theta}^G$ is contained in $\Psi^G$. 
\end{lemma}

\noindent{\bf Proof:} We first show that image of $\tilde{\theta}^G$ satisfies M1) and M2). To show M1) we see that the following diagram commutes because the basis of $\Z_p[H_1]$ as a $\Z_p[H]$-module can be taken to be the same as the basis of $\Z_p[H_1^{ab}]$ as a $\Z_p[H/[H_1,H_1]]$-module. 
\[
\xymatrix{ \tilde{K}_1(\Z_p[G]) \ar[d]_{norm^G_{H_1}} \ar[rrd]^{norm^G_H} & & \\
\tilde{K}_1(\Z_p[H_1]) \ar[rr]_{norm^{H_1}_H} \ar[d] & & \tilde{K}_1(\Z_p[H]) \ar[d] \\
\tilde{K}_1(\Z_p[H_1^{ab}])  \ar[rd]_{nr_{H_1,H}} & & \tilde{K}_1(\Z_p[H^{ab}]) \ar[ld]^{\pi_{H,H_1}} \\
& \tilde{K}_1(\Z_p[H/[H_1,H_1]])}
\]
M2) is clear. We have already shown M3) is lemma \ref{easycongruence}. Since we know M3)
\[
log\Big(\frac{\tilde{\theta}^G_H(x)^{p|H|}}{u_{G,H}(\tilde{\theta}^G(x))}\Big) 
\] 
is defined for every $H \leq G$. Using the relation
\[
\beta^G_H(L(x)) = \frac{1}{p|H|} log\Big(\frac{\tilde{\theta}^G_H(x)^{p|H|}}{u_{G,H}(\tilde{\theta}^G(x))}\Big),
\] 
and the description of image of $\beta^G$ we conclude that 
\[
log\Big(\frac{\tilde{\theta}^G_H(x)^{p|H|}}{u_{G,H}(\tilde{\theta}^G(x))}\Big)  \in p|H| T_H,
\]
for every $H \in \C$. But $log$ induces an isomorphism between $1+p|H|T_H$ and $p|H|T_H$. Hence we get M4). This shows that image of $\tilde{\theta}^G$ is contained in $\Psi^G$. \qed

\begin{definition} We define the map $\mathcal{L}: \Psi^G \rightarrow \Phi^G$ by
\[
\mathcal L((x_H)) = \Big(log\Big(\frac{x_P^{p|P|}}{u_{G,P}((x_H))}\Big)\Big)_{P \leq G}.
\]
\end{definition}

\noindent We then have the following commuting diagram
\[
\xymatrix{ \tilde{K}_1(\zp[G]) \ar[d]_{\tilde{\theta}^G} \ar[r]^{L} & \zp[Conj(G)] \ar[d]^{\beta^G} \\
\Psi^G \ar[r]_{\mathcal L} & \Phi^G}
\]

\begin{lemma} The image of $\mu_{p-1}$ embeds in $\Psi^G$ diagonally and is precisely the kernel of $\mathcal L$.
\end{lemma}
\noindent{\bf Proof:} We use induction on the order of $G$. Let $z$ be a central element of $G$ of order $p$. Let $\overline{G} = G/\langle z \rangle$. We assume by induction hypothesis that the kernel of $\Psi^{\overline{G}}$ is $\mu_{p-1}$. Consider the following commutative diagram
\[
\xymatrix{ \Psi^{G} \ar[r]^{\mathcal L}  \ar[d]_{\pi}  & \Phi^G \ar[d] \\
\Psi^{\overline{G}} \ar[r]_{\overline{\mathcal L}} & \Phi^{\overline{G}} }
\]
For any $H \leq G$ we denote $\overline{H}$ to be the image of $H$ in $\overline{G}$. The kernel of $\zp[H^{ab}]^{\times} \rightarrow \zp[\overline{H}^{ab}]^{\times}$ is trivial if $z \notin H$ and it is $1+ (\overline{z}-1)\zp[H^{ab}]$ if $z \in H$ (recall that here $\overline{z}$ is the image of $z$ is $H^{ab}$). As $z$ has order $p$, we note that 
\[
\alpha_H(x_H) =1,
\]
for any $x_H \in 1+(\overline{z}-1)\zp[H^{ab}]$. If $(x_H) \in kernel(\pi)$, then
\[
\mathcal L((x_H)) = (log(x_H^{p|H|}))_{H \leq G}.
\]
Hence $\mathcal L((x_H)) = 0$ in $\Phi^G$ for $(x_H) \in ker(\pi)$ if and only if $x_H =1$ for every $H \leq G$. Hence kernel of $\mathcal L$ injects into kernel of $\overline{\mathcal L}$ which proves the lemma.  \qed \\

\begin{lemma} The map $\tilde{\theta}^G$ is injective.
\end{lemma} 
\noindent{\bf Proof:} Let $x \in \tilde{K}_1(\zp[G])$ be in the kernel of $\tilde{\theta}^G$. Then $\beta^G(L(x)) = 0$. Since $\beta^G$ is injective we deduce $L(x) = 0$, i.e. $x\in \mu_{p-1}$. As $\mu_{p-1}$ maps isomorphically onto $\mu_{p-1} \subset \Psi^G$. Hence $x=1$.
\qed \\

\begin{lemma} The map $\tilde{\theta}^G$ is surjective. 
\end{lemma}
\noindent{\bf Proof:} Let $(x_H) \in \Psi^G$. Since $\mathcal L((x_H)) \in \Phi^G$, we get a unique $a \in \zp[Conj(G)]$ such that $\beta^G(a) = \mathcal L((x_H))$. We claim that $\omega_G(a) = 1$. This is because 
\[
\omega_G(a) = \omega_{G^{ab}}(log(\frac{x_G^{p|G|}}{u_{G,G}((x_H))})) = \omega_{G^{ab}}(log(\frac{x_G^{p|G|}}{\varphi(x_G)^{|G|}})) = 1.
\]
Hence there is a $x' \in \tilde{K}_1(\zp[G])$ such that $\mathcal L(\tilde{\theta}^G(x)) = \mathcal L((x_H))$. Hence $\tilde{\theta}^G(x')^{-1}((x_H)) \in \mu_{p-1}$. Since $\mu_{p-1} \subset \tilde{K}_1(\zp[G])$ maps isomorphically onto $\mu_{p-1} \subset \Psi^G$. We may now modify $x'$ to get an $x$ such that $\tilde{\theta}^G(x) = ((x_H))$. \qed \\

\noindent Hence we get our main theorem 

\begin{theorem} The map $\tilde{\theta}^G : \tilde{K}_1(\Z_p[G]) \rightarrow \prod_{H \leq G} \tilde{K}_1(\mathbb{Z}_p[H^{ab}])$ induces an isomorphism between $\tilde{K}_1(\mathbb{Z}_p[G])$ and $\Psi^G$.
\label{maintheorem}
\end{theorem}

\bibliographystyle{plain}
\bibliography{mybib2}

\end{document}